\documentclass[12pt]{article}
\usepackage{amsfonts,amsmath,amssymb,amsthm}

\newcommand{\de}{\delta}
\newcommand{\G}{\Gamma}
\newcommand{\g}{\gamma}
\newcommand{\te}{\theta }
\newcommand{\la}{\lambda }
\newcommand{\si}{\sigma}
\newcommand{\f}{\varphi}
\newcommand{\z}{\zeta}
\newcommand{\ep}{\varepsilon}

\newcommand{\tc}{{\mathcal T}}

\newcommand{\bc}{{\mathcal B}}

\newcommand{\exn}{{\bf E}\,}
\newcommand{\exnp}{{\bf E}}
\newcommand{\pr}{{\bf P}\,}
\newcommand{\prp}{{\bf P}}
\newcommand{\R}{\mathbb{R}}
\newcommand{\Nb}{\mathbb{N}}
\newcommand{\Zb}{\mathbb{Z}}
\newcommand{\pto}{\xrightarrow{p}}
\newcommand{\dto}{\xrightarrow{d}}
\newcommand{\deq}{\stackrel{d}{=}}

\theoremstyle{plain}
\newtheorem{theo}{Theorem}%[section]
\newtheorem{lemo}{Lemma}%[section]
\newtheorem{coro}{Corollary}%[section]

\theoremstyle{remark}
\newtheorem{rema}{Remark}%[section]

\begin{document}

\title{On the asymptotic behaviour of random recursive trees in random
environment}

\author{K.A. Borovkov\thanks{Department of Mathematics and Statistics, The University of Melbourne, Parkville 3010, Australia
}, \ \ V.A. Vatutin\thanks{Steklov Mathematical Institute RAS, Gubkin St. 8,
119991 Moscow, Russia}}

\maketitle

\begin{abstract}
We consider growing random recursive trees in random environment, in which at
each step a new vertex is attached (by an edge of a random length) to an
existing tree vertex according to a probability distribution that assigns the
tree vertices masses proportional to their random weights. The main aim of the
paper is to study the asymptotic behaviour of the distance from the newly
inserted vertex to the tree's root and that of the mean numbers of outgoing
vertices as the number of steps tends to infinity. Most of the results are
obtained under the assumption that the random weights have a product form with
independent identically distributed factors.\\


\end{abstract}

\section{Introduction}
\label{Sec1}

We consider the following random recursive tree model. A recursive tree is
constructed incrementally, by attaching a new vertex to a randomly chosen
existing tree vertex at each step.  Initially, the tree consists of a single
vertex $v(0)$ that has weight $w(0)=1$ and label $0$. At the first step, a new
vertex $v(1)$ is added to the tree as a child of the initial vertex. It is
labelled by $1,$ and a random weight $w(1)> 0$ and a random length $Y(1)\ge 0$
are assigned to the vertex and to the edge connecting the vertices $v(0)$ and
$v (1)$, respectively. It is assumed that the edge is directed from $v(0) $ to
$v(1)$. At step $j>1$, given all the weights $w(0),w(1),\dots ,w(j-1),$ first a
node $v( j^*)$ is chosen at random from the nodes $v(0),v(1),\dots ,v(j-1)$
according to the distribution with probabilities proportional to the nodes'
weights, and then a new vertex $v(j)$ is added to the tree as a child of the
node $v(j^*)$. The new vertex has label $j$, and a random weight $w(j)> 0$ and
a random length $Y(j)\ge 0$ are assigned to the vertex $v(j)$ and to the edge,
connecting the vertices $v(j^*)$ and $v(j),$ respectively. As at the initial
step (where, for consistency, we will put $1^*=0$), the edge is directed from
$v(j^*)$ to its child vertex $v(j).$ We assume that $ \{Y(j)\}_{j\ge 1}$ is a
sequence of independent random variables (r.v.'s) which is independent of the
sequence of the (generally speaking, random) weights $\{ w(j)\}_{j\ge 0}.$
Interpreting the sequence of weights as a ``random environment" in which our
recursive tree is growing, and appealing to an analogy with random walks and
branching processes in random environments, it is not unnatural to refer to
such a model as a random recursive tree in random environment.

Let
$$
D_{0}:=0,\qquad D_{n}:=D_{n^*}+Y (n), \qquad n\ge 1,
$$
be the distance from the vertex $v(n)$ to the root (i.e. the sum of the lengths
of the edges connecting $v(n)$ with $v(0)$). In this paper, we study the
asymptotic (as $n\to\infty$) behavior of $D_{n}$ under various assumptions on
the random weights $w(j)$ and lengths $Y(j)$, and also that of the mean values
of the outgoing degrees
\begin{equation}
N_n (j) := \sum_{k=j+1}^n I\{v(k^*) = v(j)\}, \qquad j\le n,
\label{out_deg}
\end{equation}
where $I\{A\}$ is the indicator of the event $A$.

Observe that if $w(j)\equiv Y(j)\equiv 1$ for all $j$, then we get the standard
random recursive tree (\cite{Dev88}; see also~\cite{Mah92}). If $w(j)=a^{j},$
$j\ge 0,$ where $a>0$ is a constant and $Y(j),$ $j\ge 1,$ are r.v.'s  whose
distributions satisfy certain mild conditions, we get the recursive tree
considered in~\cite{BoM05} (in fact, the model in~\cite{BoM05} assumed that, at
each step, a fixed number $k\ge 1$ of children are attached to one of the
existing tree vertices, and also that $Y(j)$ are vector-valued).

One should also mention here other related models where the weights of the
vertices can {\em change\/} at each step. Thus, if, after the completion of the
$k$th step of the tree construction, the weight of the vertex $v(j),$ $j\le k,$
is $w(j)=w(j,k)=1+\beta N_k (j)$, where $\beta \ge 0$ is constant and
$Y(j)\equiv 1,$ we get the linear recursive tree studied in~\cite{Pit94,BiG97}
(see also the bibliography there for further references). The case when $
w(j)=w(j,k)=1+ N_k (j)$ was considered in~\cite{BAS}; the power-tail limiting
behavior of the degree distribution for this model that had been guessed
in~\cite{BAS} was established in~\cite{BRS}.

If $w(j)=a_{1}\cdots a_{j},$ $j\ge 1,$ where  $a_{1},\dots ,a_{j}$ are
independent and identically distributed (i.i.d) r.v.'s, and $Y(j)\equiv 1,$ we
get a version of a weighted recursive tree. It is this last model and its
generalizations that will be of the main interest for us in the present paper.

From now on we assume that the weight $w(j)$ of the vertex $v(j)$ is, generally
speaking, random and, once assigned, remains unchanged forever.

Section~\ref{Sec2} of the paper is devoted to studying the asymptotic behaviour
of the distribution of $D_n$. Theorems~\ref{Tsubcrit} and~\ref{Tcritic1}
present general convergence results for the conditional distribution of $D_n$
in the cases when the random weights $w(j)$ tend to ``prescribe" new
attachments to vertices close to the root of the tree and when the new
attachments are ``more dispersed" across the tree, respectively.
Corollary~\ref{TBoM} covers the special case when $w(j)\equiv 1$. The results
of the section also show that, for any $\alpha\in (0,1],$ one can construct a
random recursive tree such that $D_n$ behaves as $n^\alpha$ as $n\to\infty$.
Theorem~\ref{Tquench1} implies that, in the case of the ``product-form" weights
$w(j)=a_{1}\cdots a_{j},$ $j\ge 1,$ with $a_j$ being non-degenerate i.i.d.
satisfying the moment conditions  $\exn \ln a_j=0$ and $\exn |\ln
a_j|^{2+\de}<\infty$ for a $\de
>0$, the limiting distribution of $D_n/\sqrt{n}$ coincides with the law of the
maximum of the Brownian motion process on a finite time interval.

Section~\ref{Sec3} deals with the expectations of the numbers of outgoing
degrees in the case of the product-form weights under the assumption that the
random walk generated by the i.i.d. sequence $\{\ln a_j\}$ satisfies Spitzer's
condition. Theorem~\ref{Tanneal1} gives the asymptotic behaviour of the
unconditional expectations $\exn N_n (k)$ as $n\to\infty$ when either $k=j$ or
$k=n-j$ for a fixed value $j\ge 0$ (in both cases it is given by a regularly
varying function of $n$). Theorem~\ref{Tanneal3} complements it by covering the
case when $\min\{j,n-j\}\to\infty$ (here the answer has the form of a product
of regularly varying functions of $j$ and $n-j$, respectively; in particular,
in the case when $\ln a_j$ has zero mean and a finite variance, one obtains
$\exn N_n (j) \sim 2\pi^{-1} (n-j)^{1/2} j ^{-1/2}$). Theorem~\ref{Texpect1}
describes, in a range of $j$-values, the asymptotic behaviour of the
distribution of the conditional expectation $\exnp_w N_n (j)$ given the
sequence of the weights $w(1), w(2),\dots.$

\section{The distribution of $D_{n}$}
\label{Sec2}

\subsection{The basic properties of $D_n$}

Let
$$
W_n:= \sum_{j=0}^{n} w(j),
 \qquad
p_{n}(j):=\frac{w(j)}{W_n},\quad j=0,1,\dots ,n.
$$
Set $f_0 (t):= 1,$ $f_{j}(t):= \exn e^{itY(j)},$ $j\ge 1,$ and put $\f_{0}(t):=
1,$
\begin{align*}
\f_{n}(t)    & := \exnp_{w}e^{itD_{n}}
 := \exn\bigl[ e^{itD_{n}}\big|\, w(1), \dots, w(n-1) \bigr],\\
\Psi _{n}(t) & := \exn\f_{n}(t) = \exn e^{itD_{n}}, \qquad n\ge 1
\end{align*}
(here and in what follows, $\exnp_{w}$ and $\prp_{w}$ denote the conditional
expectation and probability given the sequence of weights $\{w(j)\}$,
respectively).

It is easy to see that
\begin{align}
\f_{n+1}(t)
 &= \sum_{j=0}^{n} p_{n}(j)\f_{j}(t)f_{n+1}(t)  \notag \\
 &= \frac{W_{n-1}}{W_n} \sum_{j=0}^{n-1} p_{n-1}(j)\f_{j}(t)f_{n+1}(t)
  + p_{n}(n)\f_{n}(t)f_{n+1}(t)  \notag \\
 &= (1-p_{n}(n)) \frac{f_{n+1}(t)}{f_{n}(t)}\, \f_{n}(t)
  + p_{n}(n)\f_{n}(t)f_{n+1}(t)  \notag \\
 &= \bigl[ 1+(f_{n}(t)-1)p_{n}(n)\bigr] \frac{f_{n+1}(t)}{f_{n}(t)}\, \f_{n}(t)= \cdots   \notag \\
 &= f_{n+1}(t) \prod_{j=1}^{n}\bigl[ 1+(f_{j}(t)-1)p_{j}(j)\bigr] .
\label{recphi}
\end{align}

\begin{rema}
\label{rema0}
{Observe that \eqref{recphi} means in fact that, given the environment, the
r.v.  $D_{n+1}$ admits a representation in the form of a sum of independent
r.v.'s as follows:
\begin{equation}
D_{n+1} \deq I_1 Y(1) +\cdots + I_n Y(n) + Y(n+1),
\label{recphi+1}
\end{equation}
where $\{I_j\}$ is a sequence of independent (of each other and also of
$\{Y(j)\}$) random indicators with $\pr (I_j = 1)= p_j (j),$ $j\ge 1$. In the
special case when $Y(j)\equiv w(j) \equiv 1,$ this representation is equivalent
to the correspondence between the quantity $D_{n}$ and the numbers of records
in an i.i.d. sequence that was used in~\cite{Dev88} (see also Section~3.6
in~\cite{RuRa95} for a discussion of a somewhat more general situation where
the representation \eqref{recphi+1} with $Y(j) \equiv 1$ holds). Note, however,
that in~\cite{Dev88} a probabilistic argument that works in that special case
only was used to derive the representation~\eqref{recphi+1} which is actually
the main tool for studying $D_n,$ whereas our approach leads directly
to~\eqref{recphi+1} and is much more general. }
\end{rema}

From the recursive relation \eqref{recphi} one can derive a number of
interesting results on the limiting behaviour of $D_n$. Note
that~\eqref{recphi} was first derived in the case when $w(j)=a^j,$ $j\ge 0,$
$Y(j)\in \R^d,$ in~\cite{BoM05} (one can easily see that this recursive formula
and the statements of Theorems~\ref{Tsubcrit}--\ref{TBoM} below remain true in
the multivariate case as well).

In particular, the relation~\eqref{recphi} immediately implies the following
assertion, describing the limiting behaviour of the conditional (given the
weights) distribution of $D_n$ when the weight sequence $\{w(j)\}$ ``suggests"
new children to attach not too far from the tree's root.

\begin{theo}
\label{Tsubcrit}
If
$$
\sum_{j=1}^{\infty} p_{j}(j)<\infty   \qquad \text{a.s.}
$$
and the distribution of $Y(n)$ has a weak limit as $n\to\infty:$
$$
\lim_{n\to\infty} f_{n}(t)=f(t),
$$
then there exists the limit
$$
\lim_{n\to \infty} \f_{n}(t) =\f_\infty (t):=
 f(t)\prod_{j=1}^{\infty }\bigl[ 1+(f_{j}(t)-1)p_{j}(j)\bigr] \qquad \text{a.s.}
$$
\end{theo}

This result, in turn, implies that $D_n\dto D_\infty$ as $n\to\infty$, where
$D_\infty$ is a proper r.v. with the characteristic function $\exn \f_\infty
(t)$.

The next assertion refers to situations where the attachment preferences are
spread ``more uniformly" across the tree.

\begin{theo}
\label{Tcritic1}
Let the sequence of r.v.'s  $\{ Y(j)\}_{j\ge 1}$ be uniformly integrable, and
let there exist a sequence $h_n\to \infty,$ $n\to \infty,$ and a r.v. $\z$ such
that the following convergence in distribution takes place as $n\to \infty:$
\begin{equation}
\z_n:=\frac{1}{h_n}\sum_{j=1}^{n}p_{j}(j)\exn Y(j)\dto \z .
\label{odin_a}
\end{equation}
Then for any $t$
\begin{equation*}
\f_{n}\biggl( \frac{t}{h_n}\biggr)
 \dto e^{it\z }.
\end{equation*}
\end{theo}

\begin{rema}
\label{rema0+0}
{One can easily see that if, instead of~\eqref{odin_a}, one has $\z_n\to \z$
a.s. for some r.v. $\z,$ then
\begin{equation*}
\lim_{n\to \infty} \f_{n}\biggl( \frac{t}{h_n}\biggr)
 =e^{it\z }\qquad \text{a.s.}
\end{equation*}
uniformly in $t$ from any compact set.
}
\end{rema}

\begin{proof}
It is not difficult to see that, due to the uniform integrability
condition, as $n\to \infty,$
$$
f_{j}\biggl( \frac{t}{h_n}\biggr) -1
 = \exn\exp \biggl\{ \frac{itY(j)}{h_n}\biggr\} -1
 = \frac{1}{h_n}\, \bigl(it\exn Y(j)+o(1)\bigr)
$$
uniformly in $j\ge 1$ and in  $t$ from any compact set. Hence, as $p_{j}(j)\le
1,$ we have by~\eqref{recphi}
\begin{multline*}
\f_{n}\biggl( \frac{t}{h_n}\biggr)
 =
f_{n}\biggl( \frac{t}{h_n}\biggr)
 \prod_{j=1}^{n-1}\biggl[ 1+\biggl( f_{j}\biggl( \frac{t}{h_n}\biggr) -1\biggr) p_{j}(j)\biggr] \\
 =
 (1+\ep_{n}(t))\exp \biggl\{ \frac{it}{h_n}\sum_{j=1}^{n}p_{j}(j) \exn
 Y(j)\biggr\},
\end{multline*}
where $\ep_{n}(t)=o_P(1)$  as $n\to \infty $. This clearly implies the
assertion of the theorem.
\end{proof}

\begin{coro}
\label{cor1}
Under the conditions of Theorem~$\ref{Tcritic1},$
\begin{equation*}
\lim_{n\to \infty }\Psi _{n}\biggl( \frac{t}{h_n}\biggr)=\exn e^{it\z },
\end{equation*}
so that $\dfrac{D_n}{h_n}\dto \z$  as $n\to \infty $.

\end{coro}

From Theorem~\ref{Tcritic1} one can also easily deduce the following result
obtained in \cite{BoM05} (note that in the special case when $Y(j)\equiv 1$ the
result was originally established in~\cite{Dev88}).

\begin{coro}
\label{TBoM}
If~$w(j)\equiv 1,$ $j=0,1,2,\dots ,$  the family of r.v.'s $\{Y(j)\}_{j\ge 1}$
is uniformly integrable and, as $n\to \infty,$
\begin{equation*}
\frac{1}{n}\sum_{j=1}^{n}\exn Y(j) \to \mu\in \R ,
\end{equation*}
then $ \dfrac{D_n}{\ln n}\pto \mu $.
\end{coro}

\begin{proof}
In this case clearly $p_j(j) = 1/(j+1),$ and, as it was shown in Lemma~1(i)
in~\cite{BoM05}, under the above conditions
\begin{equation*}
\z_n=\frac{1}{\ln n}\sum_{j=1}^{n}\frac{1}{j+1} \, \exn Y(j)\to \mu .\qquad
\end{equation*}%
Now the assertion  of the theorem follows from
Theorem~\ref{Tcritic1}.
\end{proof}

We also get the same asymptotics for $D_n$ when the weights are random, but
remain ``on the average" the same.

\begin{coro}
\label{TBoM+0}
If~$Y(j)\equiv 1,$ $j\ge 1,$  and the sequence of random weights $\{w(j)\}$
satisfies the strong law of large numbers: as $n\to\infty$,
$$
\frac{1}{n}\sum_{j=1}^n w(j)\to a >0 \quad\text{a.s.},
$$
then $ \dfrac{D_n}{\ln n}\pto 1 $.
\end{coro}

\begin{proof}
It again suffices to apply (a slightly modified version) of Lemma~1(i) from
~\cite{BoM05} (this time to the sequences $y_n := anW_n^{-1}$, $x_n :=w(n)/a$)
and use our Theorem~\ref{Tcritic1}.
\end{proof}

\begin{rema}
\label{rema0+0+0}
{In fact, to obtain a faster than logarithmic growth rate for $D_n$ (assuming
that $Y(j)\equiv 1$), the weights $w(j)$ should grow faster that any power
function. Indeed, if, say,
$$
w(j) = j^\alpha l(j), \qquad \alpha\in\R,
$$
is a regularly varying function, then for $\alpha <-1$ one clearly has
$$
\sum_{j=1}^\infty p_j (j) <\infty
$$
(so that in this case Theorem~\ref{Tsubcrit} is applicable), whereas for
$\alpha >-1$ by Karamata's theorem $W_n \sim (\alpha +1)^{-1} n^{\alpha+1}
l(n)$, so that in this case $p_j (j) \sim 1/(\alpha+1)j$ and hence
$$
\sum_{j=1}^n p_j (j) \sim \frac{\ln n}{\alpha +1} .
$$
Thus, in the latter case $ \dfrac{D_n}{\ln n}\pto \dfrac{1}{\alpha +1} $ as
$n\to\infty$.

On the other hand, for, say,
$$
w(j) =  \alpha j^{\alpha-1} e^{j^\alpha}, \qquad \alpha\in (0,1],
$$
we get $W_n \sim e^{n^\alpha}$ and hence
$$
\sum_{j=1}^n p_j (j) \sim n^\alpha.
$$
So this example shows that, for any $\alpha\in (0,1],$ one can construct a
random recursive tree with $ \dfrac{D_n}{ n^\alpha}\pto 1 $ as $n\to\infty$.

 }
\end{rema}

\subsection{The case of the product-form random weights}

In this subsection we will construct and study recursive trees with random
vertex weights of the form $w(j) = a_1\cdots a_j,$ $j\ge 1,$ where $a_j$ are
i.i.d. r.v.'s, and unit edge lengths. As it will be clearly seen from the
proofs below, the main results will still hold true in the case of random
i.i.d. edge lengths $Y(j)\ge 0$ with a finite mean as
well~(Remark~\ref{rema1}). Thus restricting our attention to the case of unit
edge lengths leads to no loss of generality, but makes the exposition more
compact and transparent.

Denote by $\tc_{n},$ $n=0,1,2,\dots ,$ the set of all rooted recursive trees
having $n$ nonrooted vertices and unit edge lengths (that is, $\tc_{n}$
consists of the rooted trees whose root is labelled by $0$ and whose nonrooted
vertices are labelled by numbers $1,2,\dots,n$ in such a way that for any
nonrooted vertex labelled, say, by $j,$ the shortest path leading from it to
the root traverses only the vertices whose labels are less than $j$). For a
tree $t_{n}\in \tc_{n}$, let $t_{n}(j)\in \tc_{n+1}$ be the recursive tree
which is obtained from $t_{n}$ by adding a vertex labelled by $n+1$ as a child
of the vertex with the label $j\in  \{ 0,1,\dots ,n \}.$

One can describe the construction of our random recursive tree as follows.
First, we run a random walk
$$
S_{0}=0,\qquad S_{j}= \te_{1}+\cdots +\te_{j},\quad j\ge 1,
$$
where $\te_{j}\deq\te,$  $j=1,2,\dots ,n,$ are i.i.d. r.v.'s. Second,  given
$S_{j},$ $j=0,1,\dots ,n,$ we construct a (conditional) Markov chain
$T_{0},T_{1},\dots ,T_{n}$  with $T_{k}\in \tc_{k},$ $k=0,1,\dots ,n,$ by
assigning the weight
$$
w(j):=e^{-S_{j}}
$$
to the vertex labelled by $j\ge 0$ (so that $w(j)=a_1\cdots a_j,$ $j\ge 1,$ in
the notation of Section~\ref{Sec1}, with $a_j := e^{-\te_j}$ being i.i.d.
r.v.'s), so that now we have, for $r=0,1,\dots ,n,$
\begin{align}
W_r & \equiv \sum_{q=0}^{r} w(q) = \sum_{q=0}^{r}e^{-S_{q}}, \label{defb}
 \\
p_{r}(j) &\equiv \frac{e^{-S_{j}}}{W_r}
 =\frac{e^{-S_{j}}}{\sum_{q=0}^{r}e^{-S_{q}}},\quad j=0,1,\dots ,r,
\label{defprob}
\end{align}
and then letting, for any $t_{r}\in \tc_{r},$
\begin{multline*}
\prp_{w} \bigl( T_{r+1}=t_{r}(j) \big|\, T_{r}=t_{r}\bigr)\\
 \equiv \pr \bigl( T_{r+1}=t_{r}(j) \big|\, T_{r}=t_{r}; \, w(0),w(1),\dots ,w(r)\bigr)
 :=p_{r}(j),
\end{multline*}
$j=0,1,\dots ,r.$

The main result of this subsection is

\begin{theo}
\label{Tquench1}
If
\begin{equation}
\exn\te=0,
 \qquad
\si ^{2}:=\exn\te^{2}> 0,
 \qquad
\exn|\te|^{2+\de}<\infty\quad \text{for a \ $\de>0,$}
\label{spitzer}
\end{equation}
then, as $n\to \infty ,$
$$
\z_n:=\frac{1}{\sqrt{n}} \sum_{j=1}^{n}p_{j}(j)
 \dto \si_m \max_{0\le u\le 1}B(u),
 \qquad \si_m :=\si  \int_0^\infty
\frac{m(dy)}{y} <\infty,
$$
where $\{B(u)\}_{u\ge 0}$ is the standard Brownian motion process and the
measure $m$ is specified in the proof {\em (see~\eqref{7a})}.
\end{theo}

Together with Corollary~\ref{cor1}, the above assertion immediately yields the
following

\begin{coro}
\label{coro2}
Under the conditions of Theorem~$\ref{Tquench1},$
$$
\frac{D_n}{\sqrt{n}} \dto \si_m \max_{0\le u\le 1}B(u)\quad\text{as \
$n\to\infty.$}
$$
In other words, for any $x>0$
$$
\pr \bigl(D_n > \si_m \sqrt{n} x\bigr) \to 2 (1-\Phi(x)),
$$
where $\Phi$ is the standard normal distribution function.
\end{coro}

\begin{rema}
\label{rema1}
{It is obvious that the assertion of the corollary remains true in the case of
i.i.d. random edge lengths $Y(j)\ge 0$ with a finite mean, with the only
difference that $\si_m$ should be replaced in its formulation with  $\si_m\exn
Y(1).$}
\end{rema}

\begin{proof} Put
$$
L_{n} :=\min_{0\le k\le n}  S_{k}.
$$
Basing on the proof of Theorem~4.1 in~\cite{BBE97}, we will show that
\begin{equation}
\frac{1}{|L_n|} \sum_{j=1}^n p_j (j) \to \int_0^\infty \frac{m(dy)}{y} <\infty
\qquad \text{a.s.}
\label{spitzer+0}
\end{equation}
Since by the invariance principle
\begin{equation}
\frac{|L_{n}|}{\sqrt{n}}\dto \si \max_{0\le u\le 1}B(u) \qquad\text{as
 \ $n\to\infty$},
\label{spitzer+0+0}
\end{equation}
the assertion of the theorem will then immediately follow
from~\eqref{spitzer+0}.

First denote by
$$
\g_{0}:=0,\qquad \g_{j+1}:=\min \{n>\g_{j}:S_{n}<S_{\g_{j}}\}, \quad j\ge 0,
$$
the strict descending ladder epochs of the random walk $\{S_{n}\}_{n\ge 0}$.
All the r.v.'s  introduced are finite a.s. as $\{S_{n}\}_{n\ge 0}$ is recurrent
in view of~\eqref{spitzer}.

Let $\{X_n\}_{n\ge 0}$ be a Markov chain defined for $n=1,2,\dots$ by
$$
X_{n}  :=e^{\te_{n}} X_{n-1} + 1 .
$$
When $X_{0}^{x}=x>0$ is a fixed value, we will use
notation~$\{X_{n}^{x}\}_{n\ge 0}.$ Clearly,
\begin{equation}
X_{n}^{x} = x e^{S_{n}} + \sum_{q=1}^{n}e^{S_{n}-S_{q}}
 =e^{S_{n}}( x-1+W_{n}).
\label{6a}
\end{equation}

Set $\g :=\g_{1}$. Under our assumptions~\eqref{spitzer}, the expectation $\exn
S_{\g}<0$ is finite (see e.g. Corollary~10, \S~17 in~\cite{Bo76}), and the
Markov chain $\{X_{\g_{n}} \}_{n\ge 1}$ with the transition kernel
$$
M_\g (x, A) :=\pr (X_{\g }^{x}\in A),\qquad x>0,\quad A\in\bc,
$$
has a unique invariant probability measure $m_{\g }$ (see e.g. Lemma~5.49 in
\cite{Eli82} and p.481 in~\cite{BBE97}):
$$
m_{\g } ( A) =\int_0^\infty m_{\g }( dx) M_\g (x, A).
$$
Moreover,  the measure $m$ defined by
\begin{equation}
m(f):=\frac{1}{\exn(-S_{\g })} \int_0^\infty
 \exn\biggl( \sum_{k=0}^{\g -1} f(X_{k}^{x})\biggr) m_{\g }(dx)
\label{7a}
\end{equation}
is an invariant measure for the Markov chain $\{ X_{n}\}_{n\ge 0}$
(see~\cite{BBE97}).

Now note that, by virtue of~\eqref{defprob} and \eqref{6a},
$$
\z_n=\frac{1}{\sqrt{n}}\sum_{j=1}^{n} p_{j}(j)= \frac{1}{\sqrt{n}}
\sum_{j=1}^{n}\frac{1}{X_{j}^{1}}\, .
$$

Let $P_{\de_y}$ be the distribution of the two-dimensional random walk
$$
Z_{n}:=(X_{n}  ,e^{S_{n}}),\quad n\ge 0
$$
(on the group of transformations $x\mapsto ax+b$ of the real line  with the
composition law $(b_{1},a_{1}) (b_{2},a_{2}) =(b_{1}+a_{1}b_{2},a_{1}a_{2})$)
when $X_0=y$. It was shown in the proof of Theorem~4.1 of~\cite{BBE97} that if
$f\in L^{1}(m)$ then
\begin{equation}
\lim_{n\to \infty } \frac{1}{|L_{n}|}\sum_{j=1}^{n}f(X_{j})
 =\int_0^\infty f(y)\, m (dy)\qquad  P_{m_{\g }}\text{-a.s.,}
\label{conv1}
\end{equation}
where
$$
P_{m_{\g }} :=\int_0^\infty P_{\de_y} \, m_{\g } (dy)
$$
is the law of the two-dimensional random walk $\{Z_{n}\}_{n\ge 1}$ when the
distribution of $X_{0}$ is $m_{\g }.$

Let for $N=1,2,\dots $ and $x>0$
$$
g_{N}(x):=\frac{1}{x} \, I \{ N^{-1}\le x\le N \} \le \frac{1}{x}=:g(x).
$$
Clearly, for all $x>0$
\begin{equation}
g_{N}(x)\nearrow g(x) \quad\text{  as \ $N\to \infty ,$ }
\label{conas}
\end{equation}
and $g_{N}(x)\in L^{1}(m)$ for each $N=1,2,\dots .$ Therefore by~\eqref{conv1}
\begin{equation}
\lim_{n\to \infty }\frac{1}{|L_{n}|}\sum_{j=1}^{n}g_{N}(X_{j})
 =\int_0^\infty g_{N} (y)\, m (dy) \qquad P_{m_{\g }}\text{-a.s.}
\label{cutlim}
\end{equation}
On the other hand, for each $N\ge 1$ and any $x> 0$
\begin{align}
\frac{1}{|L_{n}|}&\sum_{j=1}^{n} g_N (X_{j}^{x})
  \le
\frac{1}{|L_{n}|}\sum_{j=1}^{n} g(X_{j}^{x})
 =
\frac{1}{|L_{n}|} \sum_{j=1}^{n}\frac{e^{-S_{j}}}{x-1+W_{j}}
\notag \\
 &\le
\frac{1}{|L_{n}|}\sum_{j=1}^{n}\int_{x-1+W_{j-1}}^{x-1+W_{j}} \frac{dy}{y}
 =
\frac{1}{|L_{n}|} \int_{x-1+W_0}^{x-1+W_{n}} \frac{dy}{y}
\notag \\
 &=\frac{1}{|L_{n}|} \, \bigl[ \ln (x-1+W_{n}) - \ln x\bigr]
 \le \frac{1}{|L_{n}|}\, \bigl[ \ln \bigl(x + ne^{|L_{n}|}\bigr) - \ln x\bigr]
\notag \\
 &\le  \frac{1}{|L_{n}|} \biggl[ \ln  ne^{|L_{n}|} +\frac{x}{ne^{|L_{n}|}} - \ln x\biggr]
 = 1 + \frac{1}{|L_{n}|}  \bigl[ \ln   n + O(1)\bigr] \pto 1
\label{cutsup}
\end{align}
as $n\to \infty$ by the invariance principle (see e.g.~\cite{Bi68}).

Combining \eqref{cutlim} with \eqref{cutsup} shows that
$$
\sup_{N\ge 1}\int_0^\infty  g_{N}(y)\, m(dy)\le 1,
$$
which together with~\eqref{conas} yields
\begin{equation*}
\int_0^\infty  g(y)\, m(dy)\le 1.
\end{equation*}%
Therefore by~\eqref{conv1}
\begin{equation}
\lim_{n\to \infty }\frac{1}{|L_{n}|} \sum_{j=1}^{n} g(X_{j})
 =\int_0^\infty  g(y)\, m(dy) = \int_0^\infty \frac{dm(y)}{y}\qquad
 P_{m_{\g }}\text{-a.s.}
\label{conL1}
\end{equation}
To see that this convergence holds for all starting points $x>0$, it suffices
to observe that $g(z)$ is monotone in $z>0$ and $X_{j}^{x_{1}} >
X_{j}^{x_{2}},$ $j\ge 1,$ for $x_{1}>x_{2}>0.$

This, in view of \eqref{spitzer+0} and \eqref{spitzer+0+0}, completes the proof
of Theorem~\ref{Tquench1}.

\end{proof}

\section{The expectations of the outdegrees of vertices}
\label{Sec3}

Let $N_n (j)$ be the outdegree of the vertex $v(j),\,j=0,1,\dots ,n,$ in
$\tc_n,$ i.e. the number of the edges coming out of $v(j)$ in a tree having $n$
nonrooted vertices. Clearly, the r.v. $N_n (j)$  admits the
representation~\eqref{out_deg}, and therefore
\begin{multline}
\exnp_{w} N_n (j)
 = \exn\bigl[ N_n (j) \big| \, w(1),  \dots, w(n-1)  \bigr]   \\
 = \sum_{k=j+1}^{n}\exnp_{w} I\bigl\{ v(k^*)=v(j)\bigr\}
 =\sum_{k=j+1}^{n} p_{k-1}(j)
 = e^{-S_{j}}\sum_{k=j}^{n-1} W_k^{-1}
\label{degree1}
\end{multline}
and
\begin{equation}
\exn N_n (j)=\sum_{k=j}^{n-1} \exn e^{-S_{j}} W_k^{-1}.
\label{degree2}
\end{equation}

Our aim in this section is to investigate the asymptotic (as $n\to \infty $)
behavior of the expectations $\exn N_n (j)$ and that of the distributions of
the r.v.'s $\exnp_{w} N_n (j)$ in different ranges of the parameter $j$ values.

\subsection{The asymptotic behavior of $\exn N_n (j)$}

In this section we impose weaker restrictions (compared to the
conditions~\eqref{spitzer} used in Section~\ref{Sec2}) on the random walk
$S_{n}=\te_{1}+\dots +\te_{n},$ $n\ge 1, $ where $\te_{j}\overset{d}{=}\te$ are
i.i.d. r.v.'s. Namely, we only assume that Spitzer's condition holds:

\medskip

{\em There exists a $\rho \in (0,1)$ such that
\begin{equation}
\frac{1}{n}\sum_{k=1}^{n}\pr (S_{k}>0)\to \rho \quad \text{ as } \ n\to \infty.
\label{spit}
\end{equation}
}

It is known~\cite{Don95} that this condition is equivalent to Doney's condition
\begin{equation}
\pr (S_{n}>0)\to \rho \quad \text{ as } \ n\to \infty
\label{doney}
\end{equation}
(for a further discussion of the condition~\eqref{spit}, see e.g. Section~8.9
in~\cite{BGT}).

We will need a number of auxiliary results concerning the random walk
$\{S_{n}\}_{n\ge 0}.$ Let
$$
\G_{0}:=0,\qquad \G_{j+1}:=\inf \{n>\G_{j}: \, S_{n}>S_{\G_{j}}\},\quad j\ge 0,
$$
be the strict ascending ladder epochs of the random walk $\{S_{n}\}_{n\ge 0}$.
Recall that $0=\g_0 <\g_{1} <\g_2<\dots$ denote the strict descending ladder
epochs in the walk. Introduce the two renewal functions
\begin{align*}
U(x) &:= 1+\sum_{j=1}^{\infty} \pr (S_{\G_{j}}<x),
 & x>0;\quad U( 0) =1, \qquad U(x)=0, \quad x<0, \\
V(x) &:= \sum_{j=0}^{\infty }\pr (S_{\g_{j}}\ge -x),
 & x>0;\quad V(0) =1, \qquad V( x) =0,\quad x<0,
\end{align*}
and set
$$
{M}_{n} := \max_{0\le k\le n} S_{k},
 \qquad
\widetilde{M}_{n} := \max_{1\le k\le n} S_{k}.
$$
It is known (see e.g. Lemma~1 in~\cite{HI} and Lemma~1 in~\cite{VD3}) that
under the condition~\eqref{spit}
\begin{equation}
\exn U(-\te)I \{ -\te>0\} =e^{-\phi },
 \quad
\exn U(x-\te)I \{ x-\te>0\} =U(x),
 \quad x>0,
\label{ren1}
\end{equation}
where
\begin{equation}
\phi :=\sum_{j=1}^{\infty } \frac{1}{j}\, \pr (S_{j}=0) <\infty,
\label{17a}
\end{equation}
and
\begin{equation}
\exn V(x+\te)=V(x),\quad x\ge 0.
\label{ren2}
\end{equation}
By means of $V(x)$ and $U(x)$ one can specify two sequences of probability
measures $\{ \prp_{n}^{-}\}_{n\ge 1}$ and $\{\prp_{n}^{+}\}_{n\ge 1}$ on the
$\si$-algebras $\{\Sigma_{n}:=\si (S_{1},\dots ,S_{n})\}_{n\ge 1},$
respectively, with the corresponding expectations $\{ \exnp_{n}^{-}\}_{n\ge 1}$
and $\{ \exnp_{n}^{+}\}_{n\ge 1},$ by setting for each bounded measurable
function $\psi_{n}(x_{1},\dots ,x_{n})$
\begin{equation}
\exnp_{n}^{-} \bigl[ \psi _{n} ( S_{1},\dots ,S_{n} ) \bigr]
 :=e^{\phi }\exn  \bigl[ \psi _{n} ( S_{1},\dots ,S_{n} ) U ( -S_{n} )
  I\{ \widetilde{M}_{n}<0\} \bigr]
\label{me-}
\end{equation}
and
\begin{equation}
\exn _{n}^{+}\bigl[ \psi _{n} ( S_{1},\dots ,S_{n}) \bigr]
 := \exn \bigl[ \psi _{n} ( S_{1},\dots ,S_{n}) V(S_{n})I\{ L_{n}\ge 0\} \bigr] .
\label{me+}
\end{equation}
It is easy to verify that~\eqref{ren1} and~\eqref{ren2} imply that each of the
sequences $\{\prp_{n}^{\pm}\}_{n\ge 1}$ is consistent, and therefore by
Kolmogorov's extension theorem there exist measures $\prp^{-}$ and $\prp^{+}$
on the $\si$-algebra $\si (S_{1},S_{2},\dots )$ such that their restrictions
$\prp^{\pm }|_{\Sigma _{n}}$ to $\Sigma_{n}$ coincide with $\prp_{n}^{\pm
},n=1,2,\dots .$

It is known (see Lemma~2.7 in~\cite{AGKV}) that, under the
condition~\eqref{spit},
\begin{equation}
\eta_{1} := \sum_{k=1}^{\infty } e^{S_{k}}<\infty \qquad \prp^{-}\text{-a.s.},
 \qquad
\eta_{2} := \sum_{k=0}^{\infty } e^{-S_{k}}<\infty \qquad \prp^{+}\text{-a.s.}
\label{converg}
\end{equation}
Finally, it is not difficult to deduce from Lemma~3 in~\cite{VD3} that if we
put
$$
H_{n}^{-}(x) := \pr \biggl( \sum_{k=1}^{n}e^{S_{k}} \le x \, \Big|
 \, \widetilde{M}_{n}<0 \biggr),
 \qquad
H_{n}^{+}(x) := \pr \biggl( \sum_{k=0}^{n}e^{-S_{k}}\le x \, \Big|
 \, L_{n}\ge 0\biggr),
$$
and
$$
H^{-}(x) := \prp^{-}(\eta _{1}<x),
 \qquad
H^{+}(x) := \prp^{+}(\eta _{2}<x),
$$
then under the condition~\eqref{spit}
\begin{equation}
H_{n}^{\pm }(x)\Rightarrow H^{\pm }(x)\quad\text{ as }\  n\to \infty ,
\label{convlim}
\end{equation}
where the symbol $\Rightarrow $ denotes convergence at all continuity points of
the limiting function.

In what follows we will often use the following statement (see e.g. Lemma~2.1
in~\cite{AGKV}, Theorem~8.9.12 in~\cite{BGT}, and Lemma~2 in~\cite{VD3}).

Let
$$
\la_{n} ( x) :=\pr ( L_{n}\ge -x) ,
 \qquad
\widetilde{\mu}_{n} ( x) :=\pr ( \widetilde{M}_{n}<x) ,\qquad x\ge 0.
$$

\begin{lemo}
\label{Lrenewal1}
Under Sptizer's condition~\eqref{spit} there exist slowly varying at infinity
functions $l_{1}(n)$ and $l_{2}(n)$, related by $l_{1}(n)l_{2}(n)\sim \pi
^{-1}\sin \pi \rho,$ $n\to \infty,$ such that
\begin{equation}
 \pr ( L_{n}\ge 0) \sim n^{\rho -1} l_{1}( n) ,\quad
 \pr (\widetilde{M}_{n}<0) \sim n^{-\rho} l_{2}( n) \quad\text{as \ $n\to \infty$.}
\label{asmaxs}
\end{equation}
Moreover, there are absolute constants $C_{1}>0,$ $C_{2}>0$ such that for all
$n\ge 1$ and $x\ge 0$
\begin{equation}
 \la_{n}( x) \le C_{1} V( x)\, \pr ( L_{n}\ge 0) ,\qquad
 \widetilde{\mu}_{n} (x) \le C_{2}U( x)\, \pr ( \widetilde{M}_{n}<0) .
\label{upmaxim}
\end{equation}
\end{lemo}

In~\eqref{asmaxs} and in the rest of the paper, notation $a_n\sim b_n$ means
that $a_n/b_n\to 1$ as $n\to\infty.$

Let $\{ S_{n}^{-}\}_{n\ge 0}$ and $\{ S_{n}^{+}\}_{n\ge 0}$ be two independent
copies of $\{ S_{n}\}_{n\ge 0}$, and let
$$
L_{n}^{+} := \min_{0\le r\le n} S_{r}^{+},
 \qquad
\widetilde{M}_{n}^{-} := \max_{1\le l\le n} S_{l}^-.
$$
Introduce the probability distributions
$$
 \prp_{-,+} := \prp^{-} \times \prp^{+},\quad
 \prp_{\cdot ,+} := \pr \times \prp^{+},\quad
 \prp_{-,\cdot } := \pr ^{-}\times \pr
$$
on the sample space $\R^{\infty }\times \R^{\infty }$ of the pair $\bigl(\{
S_{n}^{-}\}_{n\ge 0},\, \{ S_{n}^{+}\}_{n\ge 0}\bigr),$ where $\pr$ is the
distribution of the original sequence $\{ S_{n}\}_{n\ge 0}$ and the measures
$\prp^{\pm }$ are specified by~\eqref{me-}, \eqref{me+}, and let $\exnp_{-,+},$
$\exnp_{\cdot ,+},$ and $\exnp_{-,\cdot }$ be the expectation operators under
the respective measures.

We will call an array of r.v.'s $\{G_{l,r};\, l,r\in \Nb\}$ {\em adapted\/} if,
for any pair of indices $l,r\in\Nb,$ the r.v. $G_{l,r}$ is measurable with
respect to the $\si$-algebra $\si (S_{1}^{-},\dots ,S_{l}^{-})\otimes \si
(S_{1}^{+},\dots ,S_{r}^{+})$. The following result is contained in Lemma~3
in~\cite{VD3}.

\begin{lemo}
\label{Lfran}
Let Spitzer's condition \eqref{spit} hold, and let $\{G_{l,r}; \,l,r\in \Nb\}$
be an adapted array of uniformly bounded r.v.'s. If the following limit exists:
$$
\lim_{l,r\to \infty }G_{l,r}=:G \qquad\text{$\prp_{-,+}$-a.s.,}
$$
then
\begin{equation}
\lim_{l,r\to \infty} \exn \bigl[G_{l,r} \,|\,\widetilde{M}_{l}^{-}< 0 ,\,
L_{r}^{+}\ge 0\bigr] =\exnp_{-,+} G.
\label{g22}
\end{equation}
\end{lemo}

The next statement is a slight modification of Lemma~2.5 in~\cite{AGKV} and can
be proved by the same arguments as used there.

\begin{lemo}
\label{Lfran2}
Let Spitzer's condition \eqref{spit} hold, and let $\{G_{l,r},\,l,r\in \Nb\}$
be an adapted array of uniformly bounded r.v.'s. If the following limit exists:
$$
\lim_{r\to \infty } G_{l,r}I\{\widetilde{M}_{l}^{-}<0\}
 =:G_{l}^{+}I\{\widetilde{M}_{l}^{-}<0\} \qquad\text{$\prp_{\cdot ,+}$-a.s.,}
$$
then
$$
\lim_{r\to \infty} \exn \bigl[ G_{l,r}I\{\widetilde{M}_{l}^{-}<0\} |\,
 L_{r}^{+}\ge 0\bigr]
 =\exnp_{\cdot ,+} G_{l}^{+}I\{\widetilde{M}_{l}^{-}<0\} ,
$$
and if
$$
\lim_{l\to \infty }G_{l,r}I\{L_{r}^{+}\ge 0\}=:G_{r}^{-}I\{L_{r}^{+}\ge 0\}
\qquad\text{$\prp_{-,\cdot }$-a.s.,}
$$
then
$$
\lim_{l\to \infty } \exn \bigl[G_{l,r}I\{L_{r}^{+}\ge
0\}|\,\widetilde{M}_{l}^{-}<0\bigr]
 = \exnp_{-,\cdot }  G_{r}^{-}I\{L_{r}^{+}\ge 0\} .
$$
\end{lemo}

The following result was proved in Lemma 2.2 of~\cite{AGKV}. Denote by
$$
\tau (n):= \min \bigl\{ k\ge 0: \, S_{k}\le S_{l},\, l\in [0,n]\bigr\}
$$
the left-most point at which the random walk $\{S_n\}$ attains its minimum
value on the time interval $[0,n].$

\begin{lemo}
\label{Lconv1}
Let Spitzer's condition~\eqref{spit} hold, and let $u(x)\ge 0,$ $x\ge 0,$ be a
nonincreasing function such that $\int_{0}^{\infty }u(x)dx<\infty .$ Then, for
every $\ep >0$, there exists an integer $J$ such that for all $n\ge J$
\begin{equation*}
\sum_{p=J}^{n} \exn \bigl[ u(-S_{p});\, \tau (p) = p \bigr]
  \pr (L_{n-p}\ge 0) \le \ep \pr (L_{n}\ge 0).
\end{equation*}
\end{lemo}

Introduce the r.v.'s
\begin{equation}
G_{r}^{+}(j) := \frac{e^{-S_{j-r}^{+}} I\{j\ge r\}+e^{S_{r-j}^{-}}I\{j<r\}}{
\sum_{p=1}^{r}e^{S_{p}^{-}}+\eta _{2}^{+}}
\label{27a}
\end{equation}
and
$$
G_{r}^{-}(j):=\frac{e^{S_{j-r}^{-}}I\{j>r\}+e^{-S_{r-j}^{+}}I\{j\le r\}}{%
\eta _{1}^{-}+\sum_{p=0}^{r}e^{-S_{p}^{+}}} \, ,
$$
where $\eta _{1}^{-}$ and $\eta_{2}^{+}$ are defined as in~\eqref{converg}, but
for the random walks $\{ S_{n}^{-}\}_{n\ge 0}$ and $\{ S_{n}^{+}\}_{n\ge 0},$
respectively. Note that $0< G_r^\pm \le 1$ and, in view of~\eqref{converg},
$G_{r}^{+}(j)$ and $G_{r}^{-}(j)$ are a.s. positive under the measures
$\prp_{\cdot ,+}$ and $\prp_{-,\cdot },$ respectively. Set
\begin{equation}
\widetilde{L}_{n}^{+} :=\min_{1\le p\le n}S_{p}^{+}
\label{lwave}
\end{equation}
and put
$$
c_{j} := \sum_{l=0}^\infty \exnp_{\cdot ,+}
G_{l}^{+}(j)I\{\widetilde{M}_{l}^{-}<0\},
 \qquad
d_{j} := \sum_{q=1}^{j}\sum_{r=0}^{\infty }\exnp_{-,\cdot }
G_{r}^{-}(q)I\{\widetilde{L}_{r}^{+}>0\} .
$$
One can easily verify that $c_{j}$ and $d_{j}$ are finite for any $j=0,1,\dots
.$ Thus,
\begin{multline*}
c_{j}
 \le j+1+\sum_{l=j+1}^{\infty} \exnp_{\cdot ,+}
 e^{S_{l-j}^{-}}I\{\widetilde{M}_{l}^{-}<0\}
 = j+1+\sum_{p=1}^{\infty} \exn e^{S_{p}}I\{\widetilde{M}_{p}<0\} \\
 = j+1+\sum_{p=1}^{\infty }\exn  e^{S_{p}}I\{S_{1}<0,\dots ,S_{p}<0\} <\infty
\end{multline*}
(see Section 17, D2 in~\cite{S}).

Now we are ready to formulate and prove the following statement.

\begin{theo}
\label{Tanneal1}
Let Spitzer's condition \eqref{spit}  hold. Then for any fixed $j\ge 0$
\begin{equation}
\lim_{n\to\infty} \frac{\exn N_n (j)}{n\pr (L_{n}\ge 0)} = \frac{c_{j}}{\rho }
\label{liman1}
\end{equation}
and
\begin{equation}
\lim_{n\to\infty}\frac{\exn N_{n}(n-j)}{\pr (\widetilde{M}_{n}<0)} = d_{j}.
\label{liman2}
\end{equation}
\end{theo}

\begin{rema}
{In view of \eqref{asmaxs}, the relations \eqref{liman1}  and \eqref{liman2}
can be rewritten as
$$
\exn N_n (j) \sim c_{j} \rho^{-1} n^{\rho } l_{1}(n),
 \quad
\exn N_{n}(n-j) \sim d_{j} n^{-\rho} l_{2}(n)\quad\text{ as } \ n\to \infty .
$$
}
\end{rema}

\begin{proof}
To prove Theorem~\ref{Tanneal1}, we have to evaluate the sum~\eqref{degree2} of
expectations of the form
\begin{equation}
 \exn e^{-S_{j}} W_k^{-1} = \sum_{l=0}^k \exn e^{-S_{j}} W_k^{-1} I \{\tau (k)=l\}.
\label{30a}
\end{equation}
The key idea both in this proof and also in that of Theorem~\ref{Tanneal3} is
quite similar to that of the Laplace method: the main contribution to the
expectation~\eqref{30a} comes from the event where $j$ is close to $\tau (k)$
(for other values of $j\le k$, the quantity $e^{-S_j}$ will typically be quite
small compared to $W_k$).

First we will show that, for each fixed $\ep >0$, there exists a $J=J(\ep )$
such that for all $j\ge 0$ and all $k\ge J+j$
\begin{equation}
\exn  e^{-S_{j}} W_k^{-1} I\{\tau (k)\ge J+j\} \le \ep \pr ( L_{k-j}\ge 0).
\label{rem1}
\end{equation}
Indeed, as $W_k \ge e^{-S_{\tau (k)}},$ we have
\begin{multline*}
\exn e^{-S_{j}}W_k^{-1} I \{ \tau (k)\ge J+j\}
 \le  \exn  e^{S_{\tau (k)}-S_{j}} I\{\tau (k)\ge J+j\}   \\
 = \sum_{p=J}^{k-j} \exn e^{S_{p+j}-S_{j}}I\{\tau (k)=p+j\}
   \le \sum_{p=J}^{k-j} \exn  e^{S_{p}}I \{\tau (k-j)=p\} \\
 = \sum_{p=J}^{k-j} \exn \bigl[e^{S_{p}} I \{\tau (p)=p\} \bigr] \pr ( L_{k-j-p}\ge 0),
\end{multline*}
and to get the desired statement it remains to apply Lemma~\ref{Lconv1} with
$u(x)=e^{-x}$.

The next step is to demonstrate that for any fixed $j\ge 0,$ $l\ge 1$
\begin{equation}
\lim_{k\to \infty} \frac{\exn  e^{-S_{j}}W_k^{-1 }I \{\tau (k)=l\}}{\pr (
L_{k}\ge 0)}
 =\exn_{\cdot ,+}  G_{l}^{+}(j)I\{\widetilde{M}_{l}^{-}<0\}.
\label{as0}
\end{equation}
But this is an easy consequence of Lemma~\ref{Lfran2}. Indeed, assume first
that $j\ge l.$ Then for the r.v.'s $G_{l,r} (j)$ defined for $r\ge j-l$ by
\begin{equation}
G_{l,k-l}(j) := \frac{e^{-S_{j-l}^{+}}}{\sum_{p=1}^{l}e^{S_{p}^{-}} +
\sum_{q=0}^{k-l}e^{-S_{q}^{+}}} \le 1, \qquad k\ge j
\label{defGR}
\end{equation}
(for $r<j-l$ one can put $G_{l,r}(j) \equiv 1$), we have
\begin{multline*}
\exn  e^{-S_{j}}W_k^{-1} I \{ \tau (k)=l\}
 = \exn  \frac{e^{S_{\tau (k)}-S_{j}}}{\sum_{p=0}^{k}e^{S_{\tau (k)}-S_{p}}} I \{ \tau (k)=l \}  \\
 = \exn G_{l,k-l}(j)I \{ \widetilde{M}_{l}^{-}<0, \, L_{k-l}^{+}\ge 0\} \\
 = \exn\bigl[ G_{l,k-l}(j)I\{\widetilde{M}_{l}^{-}<0\} \big|\, L_{k-l}^{+}\ge 0\bigr] \pr ( L_{k-l}\ge 0)
\end{multline*}
(here the second relation follows from the duality principle: we use the
``time-reversed random walk" on $[0,l]$).

It is evident that, as $k\to \infty ,$
$$
 G_{l,k-l}(j)I\{\widetilde{M}_{l}^{-}<0\}\to
 G_{l}^{+}(j)I\{\widetilde{M}_{l}^{-}<0\}\qquad\text{$\prp_{\cdot ,+}$-a.s.},
$$
and  therefore  by Lemma~\ref{Lfran2}
\begin{equation}
\lim_{k\to \infty} \exn\bigl[ G_{l,k-l}(j)I\{\widetilde{M} _{l}^{-}<0\}\big|\,
L_{k-l}^{+}\ge 0 \bigr]
 =\exnp_{\cdot ,+} G_{l}^{+}(j)I\{\widetilde{M}_{l}^{-}<0\}.
\label{asplus}
\end{equation}
On the other hand, in view of~\eqref{asmaxs}  for each fixed $l$
\begin{equation}
\lim_{k\to \infty }\frac{\pr \left( L_{k-l}\ge 0\right) }{%
\pr \left( L_{k}\ge 0\right) }=1.
\label{34a}
\end{equation}
Combining this with \eqref{asplus}  gives \eqref{as0}. The case $j<l$ can be
treated in a similar way.

Now everything is ready to complete the proof of the first part of the theorem.
It follows from \eqref{30a}, \eqref{rem1} and \eqref{as0} that, for each fixed
$j\ge 0,$
\begin{equation}
\exn  e^{-S_{j}} W_k^{-1} \sim c_{j}\pr ( L_{k}\ge 0) \quad\text{ as } \ k\to
\infty.
\label{individ1}
\end{equation}
Therefore, for a fixed $\ep >0$ there exists a $K ( \ep )<\infty$ such that for
all $K\ge K (\ep)$ and $n>K$
\begin{multline}
(1-\ep )c_{j}\sum_{k=K+1}^{n-1} \pr  ( L_{k}\ge 0)
 \le  \exn N_n (j)=\sum_{k=j+1}^{K} \exn  e^{-S_{j}}W_k^{-1} +\sum_{k=K+1}^{n-1}\exn  e^{-S_{j}}W_k^{-1}  \\
 \le ( K-j) +(1+\ep )c_{j} \sum_{k=K+1}^{n-1}\pr ( L_{k}\ge 0) .
\label{35a}
\end{multline}
By \eqref{asmaxs}  and Karamata's theorem (see e.g. Section 1.6 in~\cite{BGT})
\begin{equation}
\sum_{k=K+1}^{n-1}\pr ( L_{k}\ge 0)
  \sim \frac{n}{\rho } \,\pr  ( L_{n}\ge 0) \quad\text{ as } \ n\to \infty .
\label{individ2}
\end{equation}
This together with \eqref{35a} completes the proof of~\eqref{liman1}.

Now we will prove \eqref{liman2}. Let $\{ S_{n}^*\}_{n\ge 0} \deq \{
-S_{n}\}_{n\ge 0}$ be the ``reflected" random walk. By the duality principle,
for each fixed $q\le j$
\begin{multline}
\exn  e^{-S_{n-j}}W_{n-q}^{-1}
 =\exn \frac{e^{-S_{n-j}}}{\sum_{p=0}^{n-q} e^{-S_{n-q-p}}}  \\
 =\exn \frac{e^{S_{n-q}-S_{n-j}}}{\sum_{p=0}^{n-q}e^{S_{n-q}-S_{n-q-p}}}
 =\exn \frac{e^{-S_{j-q}^*}}{\sum_{p=0}^{n-q}e^{-S_{p}^*}}
 =\exn e^{-S_{j-q}^*}(W^*_{n-q})^{-1}
\label{36o}
\end{multline}
(with an obvious definition of $W^*_{n-q}$).

Next we set
$$
 L_{n}^* :=\min_{0\le k\le n}S_{k}^*, \qquad
 \widetilde{M}_{n}^* :=\max_{1\le k\le n}S_{k}^*
$$
and observe that, as $n\to \infty ,$
\begin{equation}
\pr ( L_{n}^*\ge 0) =\pr ( M_{n}\le 0) \sim e^{\phi } \pr (
\widetilde{M}_{n}<0).
\label{36a}
\end{equation}
Indeed, putting
$$
 \chi := \inf\{k\ge 1: \, S_k\ge 0\}, \qquad
 \widetilde{\chi} :=\inf\{k\ge 1: \, S_k > 0\},
$$
we get from the factorization identities that for $|z|<1$
$$
 1-\exn z^{\widetilde{\chi}}  = \exp\biggl\{\sum_{n=0}^\infty \frac{z^n}{n}\, \pr (S_n > 0)
 \biggr\},
 \qquad
 1-\exn z^{\chi}  = \exp\biggl\{\sum_{n=0}^\infty \frac{z^n}{n}\, \pr (S_n \ge 0)\biggr\}
$$
(see e.g. Corollary~4, \S\,16 in~\cite{Bo76}). Dividing both sides of these
identities by $1-z=e^{\ln (1-z)}$, we obtain
\begin{multline*}
\sum_{n=0}^\infty z^n \pr (M_n \le 0)
 = \sum_{n=0}^\infty z^n \pr (\widetilde{\chi} > n) = \frac{1-\exn z^{\widetilde{\chi}}}{1-z}\\
 = \exp\biggl\{ - \sum_{n=1}^\infty \frac{z^n}{n}\, \pr (S_n> 0) + \sum_{n=1}^\infty \frac{z^n}{n}\biggr\}
 = \exp\biggl\{ \sum_{n=1}^\infty \frac{z^n}{n}\, \pr (S_n \le 0)\biggr\},
\end{multline*}
and similarly
$$
\sum_{n=0}^\infty z^n \pr (\widetilde{M}_n < 0)
 = \sum_{n=0}^\infty z^n \pr ( \chi  > n)
 = \exp\biggl\{ \sum_{n=1}^\infty \frac{z^n}{n}\, \pr (S_n < 0)\biggr\}.
$$
Therefore
$$
\sum_{n=0}^\infty z^n \pr (M_n \le 0)
 = e^{\phi (z)} \sum_{n=0}^\infty z^n \pr (\widetilde{M}_n < 0),\qquad
\phi (z) := \sum_{n=1}^\infty \frac{z^n}{n}\, \pr (S_n = 0).
$$
To get \eqref{36a}, it remains to use \eqref{asmaxs} and Karamata's Tauberian
theorem (see e.g. Corollary~1.7.3 in~\cite{BGT}), noting that $\phi (z)\to \phi
$ as $z\nearrow 1$.

Now from~\eqref{individ1} and \eqref{36a} we obtain that, as $n\to \infty ,$
$$
\exn  e^{-S_{j-q}^*}(W^*_{n-q})^{-1}
  \sim c_{j-q}^*\pr \left( L_{n-q}^*\ge
0\right) \sim c_{j-q}^*e^{\phi }\pr  ( \widetilde{M}_{n}<0) ,
$$
where, with a natural definition of $\exnp_{\cdot ,+}^*$ and with
$\widetilde{L}_{r}^{+}$ defined in~\eqref{lwave}, due to the
definitions~\eqref{me-} and \eqref{me+}, one has
$$
e^{\phi } c_{j-q}^*
 = e^{\phi }\sum_{l=0}^{\infty } \exnp_{\cdot ,+}^*  G_{l}^{*+}(j-q)I\{\widetilde{M}_{l}^{*-}<0\}
 = \sum_{r=0}^{\infty } \exnp_{-,\cdot }  G_{r}^{-}(j-q) I\{\widetilde{L}_{r}^{+}>0\}.
$$
Therefore we have from~\eqref{degree2} and \eqref{36o} that, as $n\to \infty ,$
\begin{multline*}
\exn N_{n}(n-j)
 = \sum_{k=n-j}^{n-1} \exn e^{-S_{n-j}}W_k^{-1}
 = \sum_{q=1}^{j} \exn  e^{-S_{j-q}^*}(W^*_{n-q})^{-1}
 \sim \pr (L_{n}^*\ge 0) \sum_{q=1}^{j} c_{j-q}^*
  \\
 \sim \pr (L_{n}^*\ge 0) e^{-\phi } \sum_{q=1}^{j} \sum_{r=0}^{\infty }
   \exnp_{-,\cdot }  G_{r}^{-}(j-q) I\{\widetilde{L}_{r}^{+}>0\}
 \sim d_{j} \pr ( \widetilde{M}_{n}<0),
\end{multline*}
as desired. Theorem~\ref{Tanneal1} is proved.
\end{proof}

The next theorem describes the asymptotic behavior of the expectation $\exn N_n
(j) $ when $\min\{j,n-j\}\to \infty .$

\begin{theo}
\label{Tanneal3}
Let Spitzer's condition~\eqref{spit} be satisfied. Then
\begin{equation}
\lim_{j,n-j\to \infty }
 \frac{\exn N_n (j)}{(n-j)\pr (\widetilde{M}_{j}<0) \pr (L_{n-j}\ge 0)}=\frac{1}{\rho }\, .
\label{42a}
\end{equation}
\end{theo}

\begin{rema}
{In view of~\eqref{asmaxs},  the assertion of the theorem can be rewritten as
$$
\exn N_n (j) \sim \rho^{-1} j^{-\rho} l_{2}(j) ( n-j)^\rho l_{1}(n-j)
 \quad\text{as }\ j,n-j\to \infty .
$$
It follows that, for any fixed $\ep \in (0,1/2),$ we have for $t\in [\ep,
1-\ep]$
$$
\exn N_n (\lfloor nt\rfloor )\sim \frac{\sin \pi \rho}{\pi \rho}
 \biggl(\frac{1-t}{t}\biggr)^\rho \quad \text{as \ $n\to\infty$.}
$$
It is interesting to compare this with the corresponding (obvious) asymptotics
for the case when $w(j)\equiv 1$: then $\exn N_n (\lfloor nt\rfloor )\sim -\ln
t$ (of course, the functions of $t$ on the right-hand sides of the both
relations are densities on $(0,1)$).

In the case when $\exn \te = 0,$ $\exn\te^2<\infty,$ we don't even need to
bound the value $j/n$ away from $0$ and $1$: in that case, from the asymptotic
behaviour of the denominators in~\eqref{42a} (see e.g. p.94 in~\cite{Bo76}), we
get
$$
\exn N_n (j) \sim \frac{2}{\pi}\biggl({\frac{n-j}{j}}\biggr)^{1/2}
 \quad\text{as }\ j,n-j\to \infty .
$$

Note also that the assertions~\eqref{liman1}, \eqref{liman2} of
Theorem~\ref{Tanneal1} can be viewed as the ``boundary cases" of~\eqref{42a}:
there is a ``smooth transition" between these asymptotics.}
\end{rema}

We split the proof of the theorem into several steps. As we said before, the
main contribution to the expectation $\exn e^{-S_j} W_k^{-1}$ from the
sum~\eqref{degree2} comes from the event where $j$ is close to $\tau (k)$. So
first we will show that the contribution from the complementary event is
negligibly small indeed.

\begin{lemo}
\label{Lan31}
Under Spitzer's condition \eqref{spit}, for any $\ep >0$ there exists a
$J=J(\ep)<\infty$ such that for all $j\ge J$ and $k-j\ge J$
\begin{equation}
\exn\bigl[ e^{S_{\tau (k)}-S_{j}}; \, |\tau (k)-j|\ge J\bigr]
 \le \ep \pr ( \widetilde{M}_{j}<0) \pr ( L_{k-j}\ge 0) .
\label{estimtermed}
\end{equation}
\end{lemo}

\begin{proof}
Fix a $J>0$ and choose a $j\ge J$ and a $k\ge j+J$. We have
$$
\exn\bigl[ e^{S_{\tau (k)}-S_{j}}; \, |\tau (k)-j|\ge J\bigr] =R_{1}+R_{2},
$$
where
$$
 R_{1} := \sum_{t=0}^{j-J}\exn\bigl[ e^{S_{\tau (k)}-S_{j}}; \,\tau (k)=t\bigr],
 \qquad
 R_{2} := \sum_{t=j+J}^{k}\exn\bigl[ e^{S_{\tau (k)}-S_{j}}; \,\tau (k)=t\bigr].
$$
First consider $R_{2}$. For $t\ge j$ we get
\begin{align*}
\exn\bigl[ e^{S_{\tau (k)}-S_{j}};\, \tau (k)=t \bigr]
 &= \exn\Bigl[ e^{S_{t}-S_{j}};\,  \min_{0\le p\le t-1} S_{p}>S_{t}, \,
   \min_{t\le p\le k} S_{p}\ge S_{t} \Bigr]
   \\
 &= \exn\Bigl[ e^{S_{t}-S_{j}};\, \min_{0\le p\le t-1} S_{p}>S_{t}\Bigr] \pr (L_{k-t}\ge 0)
   \\
 &= \exn\Bigl[ e^{S_{t-j}};\, \max_{1\le p\le t} S_{p}<0 \Bigr] \pr (L_{k-t}\ge 0)
\end{align*}
by the duality principle. Defining for each $l\ge 0$ the shifted random walk
$$
\{ S_{p}^{(l)}:=S_{l+p}-S_{l} \}_{p\ge 0},
$$
we obtain from~\eqref{upmaxim} that
\begin{align*}
\exn\Bigl[ e^{S_{t-j}};\, \max_{1\le p\le t} S_{p}<0 \Bigr]
 &= \exn\Bigl[ e^{S_{t-j}}\pr \Bigl(\max_{1\le p\le j} S_{p}^{(t-j)}<-S_{t-j} \Big|\, S_{t-j}\Bigr)
 ; \, \max_{1\le p\le t-j} S_{p}<0\Bigr]
 \\
 &= \exn\bigl[ e^{S_{t-j}}\widetilde{\mu}_{j}(-S_{t-j});\, \widetilde{M}_{t-j}<0\bigr]
 \vphantom{\Big|}
 \\
 &\le C_{2}\, \pr ( \widetilde{M}_{j}<0) \, \exn\bigl[ e^{S_{t-j}}U(-S_{t-j});\,
\widetilde{M}_{t-j}<0\bigr] .
\end{align*}
Hence
\begin{multline*}
R_{2}  \le C_{2}\, \pr  ( \widetilde{M}_{j}<0 )
  \sum_{t=j+J}^{k} \exn\bigl[ e^{S_{t-j}}U(-S_{t-j});\, \widetilde{M}_{t-j}<0\bigr] \pr (L_{k-t}\ge 0)
    \\
 = C_{2} \, \pr ( \widetilde{M}_{j}<0)
  \sum_{p=J}^{k-j}\exn\bigl[ e^{S_{p}}U(-S_{p});\, \widetilde{M}_{p}<0\bigr] \pr (L_{k-j-p}\ge 0).
\end{multline*}
Since $U(x)$ is a renewal function, we have $U(x)=O(x)$, $x\to \infty .$ Thus,
there exists a constant $C_{3}$ such that $e^{-x}U(x)\le u(x):=C_{3}e^{-x/2}$
for all $x>0.$ Since  $\int_{0}^{\infty } u(x)\, dx<\infty,$   it follows from
Lemma~\ref{Lconv1} and the duality principle that, for every $\ep >0$, there
exists a $J_{1}=J_{1}(\ep )<\infty$ such that for all $k-j>J_{1}$
$$
\sum_{p=J_{1}}^{k-j} \exn\bigl[ e^{S_{p}} U(-S_{p});\,\widetilde{M}_{p}<0\bigr]
  \pr (L_{k-j-p}\ge 0)
 \le \frac{\ep }{2C_{2}}\pr (L_{k-j}\ge 0).
$$
Thus, for $k-j>J\ge J_{1}$,
\begin{equation}
R_{2}\le \frac{\ep }{2}\, \pr  ( \widetilde{M}_{j}<0 ) \pr (L_{k-j}\ge 0).
\label{rrem1}
\end{equation}

Now we will evaluate $R_{1}$. For $t<j$ we get
\begin{align*}
\exn\bigl[ e^{S_{\tau (k)}-S_{j}};\, \tau (k)=t\bigr]
 &= \exn\Bigl[ e^{S_{t}-S_{j}};\, \min_{0\le p\le t-1} S_{p}>S_{t}; \, \min_{t\le p\le k} S_{p}\ge S_{t}\Bigr] \\
 &= \exn\Bigl[ e^{S_{t}-S_{j}};\, \min_{t\le p\le k} S_{p}\ge S_{t} \Bigr] \pr (\widetilde{M}_{t}<0) \\
 &= \exn\Bigl[ e^{-S_{j-t}};\, \min_{0\le p\le k-t} S_{p}\ge 0\Bigr] \pr (\widetilde{M}_{t}<0),
\end{align*}
where to obtain the second relation we again used the duality principle.
Arguing as before, we see that
\begin{align*}
\exn\Bigl[ e^{-S_{j-t}}; \, & \min_{0\le p\le k-t} S_{p}\ge 0\Bigr]
   \\
 &=
 \exn\Bigl[ e^{-S_{j-t}}\pr \Bigl(\min_{0\le p\le k-j} S_{p}^{(j-t)}\ge -S_{j-t}\, \Big|\, S_{j-t}\Bigr);\,
   \min_{0\le p\le j-t} S_{p}\ge 0 \Bigr]
   \\
 &=  \exn\bigl[ e^{-S_{j-t}} \la_{k-j}(S_{j-t});\, L_{j-t}\ge 0\bigr] \vphantom{\Big|}\\
 &\le  C_{1} \, \pr  ( L_{k-j}\ge 0 ) \, \exn\bigl[ e^{-S_{j-t}}V(S_{j-t});\, L_{j-t}\ge 0\bigr] .
\end{align*}
Hence
\begin{multline*}
R_{1}  \le  C_{1}\, \pr  ( L_{k-j}\ge 0 )
   \sum_{t=0}^{j-J} \exn\bigl[ e^{-S_{j-t}} V(S_{j-t});\, L_{j-t}\ge 0\bigr] \pr (\widetilde{M}_{t}<0) \\
  = C_{1}\, \pr ( L_{k-j}\ge 0) \sum_{p=J}^{j}
   \exn\bigl[ e^{-S_{p}}V(S_{p});\, L_{p}\ge 0\bigr] \pr (\widetilde{M}_{j-p}<0).
\end{multline*}
From this bound one can deduce, using Lemma~\ref{Lconv1} and the same argument
as the one employed to evaluate $R_{2},$ that for every $\ep >0 $ there exists
a $J_{2}(\ep )<\infty$ such that for all $j>J\ge J_{2}$
\begin{equation}
R_{1} \le \frac{\ep }{2}\, \pr  ( \widetilde{M}_{j}<0 ) \pr (L_{k-j}\ge 0).
\label{rrem2}
\end{equation}
Combining~\eqref{rrem1}  with \eqref{rrem2}  and setting $J:=\max
\{J_{1},J_{2}\}$ completes the proof of Lemma~\ref{Lan31}.
\end{proof}

Next we evaluate the contributions to the expectations of interest from the
events where $\tau (k)$ is equal to a fixed number close to $j$.

\begin{lemo}
\label{Lan32}
Under Sptizer's condition~\eqref{spit}, for any fixed $r\in \Zb$
\begin{equation}
\lim_{j,k-j\to \infty }
  \frac{\exn\bigl[ e^{-S_{j}}W_k^{-1};\, \tau (k)=j+r\bigr]}{\pr (\widetilde{M}_{j}<0)
  \pr (L_{k-j}\ge 0)}
 = \exnp_{-,+}  \frac{ e^{S_{r}^{-}} I \{ r\ge 0 \} + e^{-S_{-r}^{+}}I\{ r<0\}}
 {\eta _{1}^{-}+\eta _{2}^{+}} ,
\label{exactminus}
\end{equation}%
where $\eta _{1}^{-}$ and $\eta _{2}^{+}$ are independent r.v.'s  defined as
in~\eqref{converg}, but for the independent random walks $\{S_{n}^{-}\}_{n\ge
0}$ and $\{S_{n}^{+}\}_{n\ge 0}$, respectively.
\end{lemo}

\begin{proof}
For $0\le r\le k-j$ put
$$
G_{j+r,k-j-r} :=\frac{e^{S_{r}^{-}}}{\sum_{p=1}^{j+r}e^{S_{p}^{-}}
 + \sum_{p=0}^{k-j-r}e^{-S_{p}^{+}}}.
$$
Then
\begin{align*}
\exn\bigl[ e^{-S_{j}} & W_k^{-1};\, \tau (k)=j+r\bigr]
   \\
 &= \exn\biggl[ \frac{e^{S_{j+r}-S_{j}}}{\sum_{p=0}^{k}e^{S_{j+r}-S_{p}}}; \,
  \min_{0\le p\le j+r-1} S_{p}>S_{j+r};\, \min_{j+r\le p\le k}S_{p}\ge S_{j+r}\biggr]
   \\
 &= \exn\bigl[ G_{j+r,k-j-r}; \, \widetilde{M}_{j+r}^{-}<0,\, L_{k-j-r}^{+}\ge 0\bigr]
   \\
 &= \exn\bigl[ G_{j+r,k-j-r}\, \big|\,  \widetilde{M}_{j+r}^{-}<0, \, L_{k-j-r}^{+}\ge 0\bigr]
    \pr ( \widetilde{M}_{j+r}<0) \pr (L_{k-j-r}\ge 0).
\end{align*}
Clearly, $0<G_{j+r,k-j-r}\le 1$ and
$$
\lim_{j,k-j\to \infty} G_{j+r,k-j-r}
 =\frac{e^{S_{r}^{-}}}{\eta_{1}^{-}+\eta_{2}^{+}}\qquad\text{$\prp_{-,+}$-a.s.}
$$
Hence, applying Lemma~\ref{Lfran} and recalling~\eqref{asmaxs} and the
properties of regularly varying functions (cf.~\eqref{34a}), we
get~\eqref{exactminus} for $r\ge 0$. The proof of~\eqref{exactminus}  in the
case $r<0$ is almost identical. Lemma~\ref{Lan32} is proved.
\end{proof}

\begin{proof}[Proof of Theorem~\ref{Tanneal3}]
For a fixed $\ep >0$ let $J=J(\ep )$ be such
that~\eqref{estimtermed} holds true. For $j\ge J$ and $n-j\ge J+1$ we have
from~\eqref{degree2} that
\begin{equation*}
\exn N_n (j)=R_{3}+R_{4}+R_{5},
\end{equation*}%
where
$$
R_{3}:=\sum_{k=j}^{j+J-1} \exn  e^{-S_{j}}W_k^{-1} ,
 \qquad
R_{4}:=\sum_{k=j+J}^{n-1} \exn \bigl[ e^{-S_{j}}W_k^{-1};\,
  |\tau (k)-j|<J\bigr]
$$
and
$$
R_{5}:=\sum_{k=j+J}^{n-1}\exn\bigl[ e^{-S_{j}}W_k^{-1}; \,
 |\tau (k)-j|\ge J \bigr].
$$
We evaluate the quantities $R_{i},$ $i=3,4,5,$ separately. First observe that,
in view of~\eqref{liman2} (with $n$ replaced by $k)$, there exists a constant
$C_{3}$ such that for all sufficiently large $j$
$$
R_{3}\le C_{3}J\, \pr  ( \widetilde{M}_{j}<0 ),
$$
and since
$$
(n-j) \, \pr  ( L_{n-j}\ge 0 ) \sim  (n-j)^{\rho } l_{1}(n-j)\to \infty
 \quad\text{ as }\ n-j\to \infty ,
$$
it follows that
\begin{equation}
R_{3}=o\Bigl( (n-j)\, \pr  ( \widetilde{M}_{j}<0 ) \pr  ( L_{n-j}\ge 0 ) \Bigr)
 \quad\text{ as } \ n-j\to \infty .
\label{r3}
\end{equation}

Further, using the obvious inequality $W_k \ge e^{-S_{\tau(k)}}$ and the
bound~\eqref{estimtermed} together with~\eqref{asmaxs} and Karamata's theorem,
we get for $j\ge J$ and some positive absolute constant $C_{5}$ that
\begin{multline*}
R_{5} \le \ep \pr (\widetilde{M}_j<0 ) \sum_{k=j+J}^{n-1} \pr ( L_{k-j}\ge 0)
   \\
 = \ep \pr  ( \widetilde{M}_{j}<0 )
   \sum_{p=J}^{n-j-1} \pr ( L_{p}\ge 0) \le \ep C_{5} (n-j)\,\pr( \widetilde{M}_{j}<0) \pr ( L_{n-j}\ge 0),
\end{multline*}
and therefore
$$
\frac{R_{5}}{(n-j)\,\pr (\widetilde{M}_j<0) \pr (L_{k-j}\ge 0 )}\le \ep C_{5}.
$$

Finally, set
$$
E _{J} := \exnp_{-,+} \frac{1+\sum_{r=1}^{J-1}
 \bigl( e^{S_{r}^{-}}+e^{-S_{r}^{+}}\bigr)}{\eta _{1}^{-}+\eta _{2}^{+}}\, .
$$
Using Lemma~\ref{Lan32}, the relation~\eqref{asmaxs} and the properties of
regularly varying functions, we see that, as $\min \{j,n-j\}\to \infty ,$
\begin{multline*}
R_{4}
  \sim E_{J} \pr ( \widetilde{M}_{j}<0)
    \sum_{k=j+J}^{n-1} \pr ( L_{k-j}\ge 0)
   \\
  \sim E_{J} \pr ( \widetilde{M}_{j}<0)
    \sum_{p=J}^{n-j-1}\pr ( L_{p}\ge 0)
  \sim E_{J} \pr ( \widetilde{M}_{j}<0) \rho^{-1} (n-j)\,\pr ( L_{n-j}\ge 0).
\end{multline*}
Since $\lim_{J\to\infty} E_J =1$ by the dominated convergence theorem, the
assertion of Theorem~\ref{Tanneal3} immediately follows from the above relation
for $R_4$ and the bounds for $R_3$ and $R_5$.
\end{proof}

\subsection{The asymptotic behavior of the distribution of $\exnp_{w}N_n (j)$}

Unfortunately,  our description of the asymptotic behavior of $\exnp_{w}N_n
(j)$ will be less detailed than that of $\exn N_n (j)$. We will be able to
describe the distribution of the r.v. $\exnp_{w}N_n (j)$ only for $j$ located
either to the right  or in a small left vicinity of the random epoch $\tau
(n)$.

\begin{theo}
\label{Texpect1}
Let Spitzer's condition \eqref{spit} be satisfied  and  $j=j(n)$ be an
arbitrary $($random$)$ sequence with the property that $(\tau (n)-j)_+ = o(n)$
in probability as $n\to\infty$. Then
\begin{equation}
\pr \biggl( \frac{e^{S_{j}-S_{\tau (n)}}}{n-j} \,\exn_{w} N_n (j)<x \biggr)
 \Rightarrow
\prp_{-,+} \biggl( \frac{1}{\eta _{1}^{-}+\eta _{2}^{+}} <x \biggr) ,
\label{node1}
\end{equation}%
where $\eta _{1}^{-}$ and $\eta _{2}^{+}$ are r.v.'s defined as in
\eqref{converg}, but for the independent random walks $\{S_{n}^{-}\}_{n\ge 0}$
and $\{S_{n}^{+}\}_{n\ge 0}$, respectively.
\end{theo}

\begin{proof}
Since the r.v.'s  $W_n$ (see~\eqref{defb}) are increasing in $n,$
we have from~\eqref{degree1} the following lower bound:
$$
\exnp_{w}N_n (j)
 \ge ( n-j)\, e^{-S_{j}} W_n^{-1}
 = \frac{( n-j)\, e^{S_{\tau (n)}-S_{j}}}{\sum_{k=0}^{n}e^{S_{\tau (n)}-S_{k}}}\, .
$$

Now we will derive an upper bound for $\exnp_{w}N_n (j).$ To this end observe
that, according to~\eqref{converg}, for any fixed $\ep >0$ and $\delta
>0$ there exists a $J<\infty$ such that
\begin{equation}
\prp^{+}\biggl( \sum_{k=J}^{\infty } e^{-S_{k}}>\de \biggr) \le \ep .
\label{mesurplus}
\end{equation}
Clearly, for any $j\in \lbrack \tau (n),n-1]$%
\begin{multline*}
\exnp_{w}N_n (j)
 \le  e^{S_{\tau (n)}-S_{j}}(\tau
(n)+J-j)_{+}+e^{-S_{j}}(n-j)W_{\tau (n)+J}^{-1}
   \\
 = e^{S_{\tau (n)}-S_{j}}\Biggl[ (\tau (n)+J-j)_{+}
  + (n-j)\biggl( \sum_{k=0}^{\tau (n)+J}e^{S_{\tau (n)}-S_{k}}\biggr)^{-1} \Biggr].
\end{multline*}
Hence we get
\begin{multline}
\biggl( \sum_{k=0}^{n} e^{S_{\tau (n)}-S_{k}}\biggr)^{-1}
 \le \frac{e^{S_{j}-S_{\tau (n)}}}{n-j}\, \exnp_{w}N_n (j)
    \\
 \le \frac{(\tau (n)+J-j)_{+}}{n-j}
   + \biggl( \sum_{k=0}^{\tau (n)+J} e^{S_{\tau\ (n)}-S_{k}}\biggr)^{-1}.
\label{bubul}
\end{multline}
Evidently, for $y>0$
\begin{multline}
\pr \left( \sum_{k=0}^{n}e^{S_{\tau (n)}-S_{k}}<y\right)
 =\sum_{p=0}^{n} \pr \left( \sum_{k=0}^{n}e^{S_{\tau(n)}-S_{k}}<y;\, \tau (n)=p\right)
   \\
 = \sum_{p=0}^{n}\pr \left(\sum_{l=1}^{p}e^{S_{l}^{-}}+\sum_{r=0}^{n-p}e^{-S_{r}^{+}}<y;
  \, \widetilde{M}_{p}^{-}<0, \, L_{n-p}^{+}\ge 0\right) .
\label{arcs1}
\end{multline}
Further, note that from~\eqref{converg} and~\eqref{convlim}, as $\min
\{p,n-p\}\to \infty,$
\begin{equation}
 \pr \left( \sum_{l=1}^{p} e^{S_{l}^{-}} + \sum_{r=0}^{n-p} e^{-S_{r}^{+}}<y\, \bigg|
  \, \widetilde{M}_{p}^{-}<0, \, L_{n-p}^{+}\ge 0 \right)
 \Rightarrow \prp_{-,+} \bigl(\eta _{1}^{-}+\eta _{2}^{+}<y\bigr).
\label{covlim2}
\end{equation}
If the condition~\eqref{spit} is met, then the generalized arcsine law holds
true (see e.g. Theorems~8.9.9, 8.9.5 in~\cite{BGT}):
\begin{equation}
\lim_{n\to \infty} \pr \biggl( \frac{\tau (n)}{n}\le x\biggr)
  = \frac{\sin \pi \rho }{\pi }\int_{0}^{x} t^{\rho -1} (1-t)^{-\rho}dt,
  \qquad x\in [0,1].
\label{Arksin}
\end{equation}
Thus, for any $\ep_{1}>0$ there exists a $\de_{1}\in (0,1/2)$ such that
\begin{equation}
\pr \bigl(\tau (n)\notin (n\delta_{1},n(1-\delta_{1})) \bigr)\le \ep_{1},
\label{arc3}
\end{equation}
which, combined with~\eqref{arcs1} and \eqref{covlim2}, shows that, as
$n\to\infty,$
\begin{equation}
\pr \left( \sum_{k=0}^{n}e^{S_{\tau (n)}-S_{k}}<y\right)
 \Rightarrow
\prp_{-,+} \bigl(\eta _{1}^{-}+\eta _{2}^{+}< y \bigr).
\label{below1}
\end{equation}
A similar argument combined with~\eqref{mesurplus}  shows that
\begin{equation}
\pr \left( \sum_{k=0}^{\tau (n)+J}e^{S_{\tau (n)}-S_{k}}<y\right)
 \Rightarrow
\prp_{-,+}\bigl(\eta _{1}^{-}+\eta _{2}^{+}<y\bigr)
\label{above1}
\end{equation}
as first $n\to \infty$, and than $J\to \infty.$ On the other hand, again using
\eqref{Arksin}, we conclude that, within the range $j\in \lbrack \tau
(n),n-1],$
\begin{align}
\frac{(\tau (n)+J-j)_{+}}{n-j}
 &\le  I\{ \tau (n)+J>j\} \frac{J}{n-j} I\bigl\{ \tau (n)\ge n-\sqrt{n}\bigr\}
    \notag\\
 &\qquad +\, I \{ \tau (n)+J>j \} \frac{J}{\sqrt{n}-J} I\bigl\{ \tau (n)<n-\sqrt{n}\bigr\}
    \notag\\
 &\le JI\bigl\{ \tau (n)\ge n-\sqrt{n}\bigr\} + \frac{J}{\sqrt{n}-J} \pto 0  \label{above2}
\end{align}
as first $n\to \infty,$ and than $J\to \infty.$

Using \eqref{below1} and \eqref{above1}, \eqref{above2} on the left- and
right-hand sides of \eqref{bubul}, respectively,  proves~\eqref{node1}  for
$j\in \lbrack \tau (n),n-1].$

For $\tau (n)-j>0$ one can use similar arguments. The only difference is that
in this case
$$
(\tau (n)+J-j)_{+}=\tau (n)+J-j,
$$
and for $j<\tau (n),$ varying with $n$ in such a way that $(\tau (n) -j)_+
=o(n)$, the conclusion~\eqref{above2} still holds by~\eqref{arc3}.
Theorem~\ref{Texpect1} is proved.
\end{proof}

\bigskip\noindent{\bf Acknowledgments.} Research supported by the ARC Centre of Excellence for Mathematics and
Statistics of Complex Systems. The second author was also supported by the
Russian Foundation for Basic Research (grant 05-01-00035) and by the program
``Contemporary Problems of Theoretical Mathematics" of the Russian Academy of
Sciences. He is grateful to the Department of Mathematics and Statistics of the
University of Melbourne  for its hospitality while he was visiting the
department.

%The authors are also grateful to the referee for his valuable comments that led
%to an improved exposition of the paper. }

\end{document}